\documentclass[12pt]{article}
\usepackage{amssymb}

\newtheorem{prethm}{{\bf Theorem}}

\newenvironment{thm}{\begin{prethm}{\hspace{-0.5
               em}{\bf .}}}{\end{prethm}}
\newtheorem{prelemma}{{\bf Lemma}}

\newtheorem{preex}{{\bf Example}}

\newtheorem{preprop}{{\bf Proposition}}

\newenvironment{prop}{\begin{preprop}{\hspace{-0.5em}{\bf .}}}{\end{preprop}}
\newtheorem{precor}{{\bf Corollary}}

\newenvironment{cor}{\begin{precor}{\hspace{-0.5
               em}{\bf .}}}{\end{precor}}
\newtheorem{preremark}{{\bf Remark}}

\newtheorem{preprob}{{\bf Problem}}

\newenvironment{prob}{\begin{preprob}{\hspace{-0.5
               em}{\bf.}}}{\end{preprob}}
\newtheorem{predefin}{{\bf Definition}}

\newenvironment{defin}{\begin{predefin}{\hspace{-0.5
               em}{\bf .}}}{\end{predefin}}
\newtheorem{preconj}{{\bf Conjecture}}

\newenvironment{conj}{\begin{preconj}{\hspace{-0.5
               em}{\bf .}}}{\end{preconj}}
\newtheorem{preprobb}{{\bf Problem}}

\newtheorem{prelem}{{\bf Theorem}}

\newenvironment{proof}{{\bf Proof.}\rm }{\hfill{$\Box$}}

\newtheorem{presolution}{{\bf Solution.}}

\title{\vspace{-2cm}\Large\bf\noindent First-Fit coloring of Cartesian product graphs and its defining sets}

\author{\large\bf Manouchehr Zaker\footnote{mzaker@iasbs.ac.ir}
\vspace{5mm}\\
    Department of Mathematics,\\
     Institute for Advanced Studies in Basic Sciences,\\
    Zanjan, Iran\\
  }
    \date{}
\begin{document}
\maketitle

\begin{center}
	{\large\hspace*{1.6cm}{In memory of a friend and combinatorist\newline Mojtaba Mehrabadi (1969--1998)}}
\end{center}

\vspace{0cm}\begin{abstract}
\noindent Let the vertices of a Cartesian product graph $G\Box H$ be ordered by an ordering $\sigma$. By the First-Fit coloring of $(G\Box H, \sigma)$ we mean the vertex coloring procedure which scans the vertices according to the ordering $\sigma$ and for each vertex assigns the smallest available color. Let $FF(G\Box H,\sigma)$ be the number of colors used in this coloring. By introducing the concept of descent we obtain a sufficient condition to determine whether $FF(G\Box H,\sigma)=FF(G\Box H,\tau)$, where $\sigma$ and $\tau$ are arbitrary orders. We study and obtain some bounds for $FF(G\Box H,\sigma)$, where $\sigma$ is any quasi-lexicographic ordering. The First-Fit coloring of $(G\Box H, \sigma)$ does not always yield an optimum coloring. A greedy defining set of $(G\Box H, \sigma)$ is a subset $S$ of vertices in the graph together with a suitable pre-coloring of $S$ such that by fixing the colors of $S$ the First-Fit coloring of $(G\Box H, \sigma)$ yields an optimum coloring. We show that the First-Fit coloring and greedy defining sets of $G\Box H$ with respect to any quasi-lexicographic ordering (including the known lexicographic order) are all the same. We obtain upper and lower bounds for the smallest cardinality of a greedy defining set in $G\Box H$, including some extremal results for Latin squares.
\end{abstract}




\section{Introduction}

\noindent Let $(G,\sigma)$ be a graph whose vertices are ordered by $\sigma$. The First-Fit coloring of $(G,\sigma)$ scans the vertices of $G$ according to the ordering $\sigma$, assigns the color 1 to the first vertex and at each step of the coloring, assigns the minimum available number to a vertex $v$ which has not appeared in the set of previously colored neighbors of $v$. Denote the number of colors used in First-Fit coloring of $(G,\sigma)$ by $FF(G, \sigma)$. The famous parameter Grundy number (also known as the First-Fit chromatic number) of $G$, denoted by $\Gamma(G)$ is defined as $\Gamma(G)= {\max}_{\sigma} FF(G, \sigma)$, where the maximum is taken over all orderings $\sigma$ on $V(G)$. There are many papers concerning the Grundy number and First-Fit coloring of graphs (e.g. \cite{AHL, BHLY, CGHLM, HKS, Z2}). It is clear that for any ordering $\sigma$ on the vertex set of $G$ we have $\chi(G)\leq FF(G,\sigma)\leq \Gamma(G)$. Also the inequality $\Gamma(G)\leq \Delta(G)+1$ holds, where $\Delta(G)$ is the maximum degree in $G$. But there is not any upper bound for $\Gamma(G)$ in terms of any function of $\chi(G)$. For example trees may have arbitrary large Grundy numbers. The main motivation to study the Grundy number of graphs is due to the fact that we do not know already for which orderings $\sigma$ on the vertex set of a graph $G$, the First-Fit algorithm outputs a coloring with reasonable number (with respect to $\chi(G)$) of colors. For this reason we study the worst-case behavior of the First-Fit algorithm, that is the Grundy number of graphs. Throughout the paper the complete graph on $n$ vertices is denoted by $K_n$. In this paper we study the First-Fit coloring of the Cartesian product of graphs, hence we need to present the required terminology.

\noindent Let $G=(V(G),E(G))$ and $H=(V(H),E(H))$ be any two undirected graphs without any loops or multiple edges. By the Cartesian product $G\Box H$ we mean a graph on the vertex set $V(G)\times V(H)$, where any two vertices $(u,v)$ and $(u',v')$ are adjacent if and only if either $u=u'$ and $vv'\in E(H)$ or $uu'\in E(G)$ and $v=v'$. There exists an efficient ${\mathcal{O}}(m \log n)$ algorithm such that given any graph $L$ on $n$ vertices and $m$ edges, determines whether $L$ is Cartesian product of two graphs and in this case outputs $G$ and $H$ such that $L=G\Box H$ \cite{AHI}. Let $G$ be a graph whose vertices are ordered by an ordering $\sigma$. For any two vertices $u$ and $u'$ of $G$, by $\sigma(u)<\sigma(u')$ we mean the vertex $u$ appears before $u'$ in the ordering $\sigma$. Assume that $(G,\sigma)$ and $(H, \sigma')$ are two ordered graphs. The {\it lexicographic ordering} of $G\Box H$ induced by $\sigma$ and $\sigma'$ is defined as follows. Let $(u,v)$ and $(u',v')$ are two vertices of $G\Box H$. Then in the lexicographic ordering, $(u,v)$ appears before $(u',v')$ if and only if either $\sigma(u) < \sigma(u')$ or $\sigma'(v) < \sigma'(v')$. By $(G\Box H, lex)$ we mean the graph $G\Box H$ whose vertices are ordered lexicographically. The lexicographic order is not the only ordering to be studied in this paper. But the systematic representation of $G\Box H$ is by a $|G|\times |H|$ array, where the rows are indexed by the vertices of $G$ and arranged up-down according to the ordering of $V(G)$ and the columns are indexed by the vertices of $H$ and arranged left-right according the ordering of $V(H)$.
For example while vertex coloring of $G\Box H$, the systematic way is to scan the vertices of $G\Box H$ from left to right and from up to down. Observe that this is equivalent to the scanning of $V(G\Box H)$ according to the lexicographic order. Hence the lexicographic ordering is a natural ordering of the vertices of $G\Box H$. In the following we define quasi-lexicographic orderings of $V(G\Box H)$.

\begin{defin}
Let $(G,\sigma)$ and $(H,\sigma')$ be two ordered graphs. An ordering $\tau$ on the vertex set of $G\Box H$ is called {\it quasi-lexicographic} if for any two vertices $(u,v)$ and $(u',v')$ in $G\Box H$, $\tau(u,v)<\tau(u',v')$ implies $\sigma(u)<\sigma(u')$ or $\sigma'(v)<\sigma'(v')$.
\end{defin}

\noindent Given any two graphs $(G, \sigma)$ and $(H, \sigma')$, there are many quasi-lexicographic orders corresponding to $G\Box H$ and only one of them is the lexicographic order. In fact, to determine the number of quasi-lexicographic orders is a difficult problem. For example there are 26 quasi-lexicographic orderings for the graph $K_3\Box K_3$. Arrange the vertices of $K_3\Box K_3$ lexicographically by $1, \ldots, 9$. Then $1, 4, 7, 2, 3, 5, 8, 6, 9$ is a quasi-lexicographic order. The First-Fit coloring of Cartesian product of graphs has been the research subject of many papers. Specially the Grundy number of Cartesian graph product was widely studied in the literature \cite{AHL,BHLY,CGHLM,HJ,HKS}. Unfortunately as proved by Ast\'{e} et al. \cite{AHL}, there is no upper bound for $\Gamma(G\Box G)$ in terms of any function in $\Gamma(G)$. This gives motivations to study First-Fit coloring of $G\Box H$ with certain vertex orderings on the vertex set of $G\Box H$, such as the lexicographic or quasi-lexicographic orderings.

\noindent The paper is organized as follows. In Section 2 we introduce the concept of descent and using this concept we prove that study of the First-Fit coloring with respect to quasi-lexicographic orderings is reduced to the First-Fit coloring of $(G\Box H, lex)$ (Theorem \ref{descent-free}). Then we obtain some results and bounds for $FF(G\Box H, lex)$. Section 3 devotes to greedy defining sets (to be defined later) in vertex colorings of $G\Box H$ and also Latin rectangles. In this section we first prove that study of greedy defining sets for quasi-lexicographic orderings is reduced to study of greedy defining sets in $FF(G\Box H, lex)$ (Theorem \ref{quasi-gds}). Then some upper and lower bounds are obtained in this section. In Section 4 we propose some unsolved problems and conjecture for further researches in the area of greedy defining sets.

\section{First-Fit coloring of $(G\Box H, \sigma)$}

\noindent Let $\tau$ be any quasi-lexicographic ordering for the Cartesian product of $(G, \sigma)$ and $(H, \sigma')$. In Theorem \ref{descent-free} we show that study of the First-Fit coloring of $(G\Box H, \tau)$ is reduced to study of the First-Fit coloring of $(G\Box H, lex)$. We need to introduce a key concept which we call {\it descent}.

\begin{defin}
\noindent Let $G$ and $H$ be two ordered graphs and $\tau$ be any ordering for the vertices of $G\Box H$. Let also $\mathcal{C}$ be a proper vertex coloring of $G\Box H$ using colors $1, 2, \ldots, k$. Let $x$ and $y$ be any two arbitrary colors with $1\leq x<y
\leq k$ and $v$ be an arbitrary vertex of $G\Box H$ whose color in $\mathcal{C}$ is $y$. Let $N$ be the set of neighbors of $v$ in $G\Box H$ whose color is $x$. Then $\{v\}\cup N$ is said to be a $(\mathcal{C},\tau)$-descent for $G\Box H$ if any vertex $u\in N$ satisfies $\tau(u)>\tau(v)$. In case that $N$ is the empty set then $\{v\}$ is a descent.
\end{defin}

\noindent Given $(G\Box H, \tau, \mathcal{C})$, we say $\mathcal{C}$ is descent-free if there exists no $(\mathcal{C},\tau)$-descent in $G\Box H$. The following theorem shows the relation between descent-free colorings and First-Fit coloring of $G\Box H$.

\begin{thm}
Let $\mathcal{C}$ be a descent-free coloring of $(G\Box H, \tau)$, where $\tau$ is an arbitrary ordering. Let $\mathcal{D}$ be the coloring obtained by the First-Fit coloring of $(G\Box H, \tau)$. Then $\mathcal{C}=\mathcal{D}$. In particular the number of colors used in both colorings are the same.\label{descent}
\end{thm}

\noindent\begin{proof}
Assume that $v_1, v_2, \ldots, v_n$ is an ordering of the vertices of $G\Box H$ such that $\tau(v_1)< \tau(v_2) < \ldots < \tau(v_n)$. For any $i$, denote the color of $v_i$ in $\mathcal{C}$ (resp. $\mathcal{D}$) by $\mathcal{C}(v_i)$ (resp. $\mathcal{D}(v_i)$). By the definition of First-Fit coloring, $\mathcal{D}(v_1)=1$. We claim that $\mathcal{C}(v_1)=1$. Otherwise let $\mathcal{C}(v_1)=i>1$. Let $N$ be the set of neighbors of $v_1$ having color 1 in $\mathcal{C}$. If $N= \emptyset$ then $\{v_1\}$ is descent. If $N\not= \emptyset$ then by our choice of $v_1$, any $u\in N$ satisfies $\tau(u)>\tau(v_1)$. In this case $\{v_1\}\cup N$ is descent. Hence $\mathcal{C}(v_1)=1$. Assume that $v_1, \ldots, v_k$ are so that $\mathcal{C}(v_i)=\mathcal{D}(v_i)$ for any $i\in \{1, \ldots k\}$. We prove that $\mathcal{C}(v_{k+1})=\mathcal{D}(v_{k+1})$. Set $\mathcal{C}(v_{k+1})=j$. Since $\mathcal{C}$ is descent-free then for any $i\in \{1, \ldots, j-1\}$ there exists a neighbor of $v_{k+1}$ say $v_{n_i}$ such that $\mathcal{C}(v_{n_i})=i$ and $\tau(v_{n_i})<\tau(v_{k+1})$. If follows that for any $i\in \{1, \ldots, j-1\}$, $v_{n_i}\in \{v_1, \ldots, v_k\}$. Hence $\mathcal{D}(v_{n_i})=\mathcal{C}(v_{n_i})=i$. It implies that $j$ is the first available color for coloring the vertex $v_{k+1}$ by the First-Fit procedure. In other words, $\mathcal{D}(v_{k+1})=j=\mathcal{C}(v_{k+1})$. It follows that $\mathcal{C}$ and $\mathcal{D}$ are the same colorings.
\end{proof}

\noindent In the following we apply Theorem \ref{descent} for quasi-lexicographic orderings.

\begin{thm}
Let $\mathcal{C}$ (resp. $\mathcal{C}'$) be the coloring of $G\Box H$ obtained by the First-Fit coloring of $(G\Box H, lex)$ (resp. $(G\Box H, \tau)$), where $\tau$ is any quasi-lexicographic ordering of $G\Box H$. Then $\mathcal{C}=\mathcal{C}'$ and $FF(G\Box H, lex)=FF(G\Box H, \tau)$.\label{descent-free}
\end{thm}

\noindent\begin{proof}
Since $\mathcal{C}$ is obtained from the First-Fit coloring of $(G\Box H, lex)$, then there is not any descent in $(G\Box H, lex, \mathcal{C})$. We show that $\mathcal{C}$ is descent-free in $(G\Box H, \tau, \mathcal{C})$. Assume on the contrary that for some $\alpha \in V(G\Box H)$, $\{\alpha\}\cup N$ is a $(\mathcal{C}, \tau)$-descent. Denote the color of any vertex $w$ of $G\Box H$ in $\mathcal{C}$ by $\mathcal{C}(w)$. Let $\mathcal{C}(\alpha)=y$ and $x$ be the color of any vertex in $N$. Recall the definition of descent and obtain $\tau(\beta) > \tau(\alpha)$ for any vertex $\beta \in N$. In order to obtain a contradiction, we show that $\{\alpha\}\cup N$ is a descent in $(\mathcal{C}, lex)$. For this purpose it is enough to show that for any $\beta \in N$ we have $lex(\beta) > lex(\alpha)$, where by $lex(\beta)$ we mean the order of the vertex $\beta$ in the lexicographic ordering. Let $\alpha=(u,v)$ and $\beta=(u',v')$ be any vertex of $N$, where $u, u'\in G$ and $v, v'\in H$. Since $\tau$ is quasi-lexicographic and $\tau(u',v') > \tau(u,v)$, then either $u'>u$ in $G$ or $v'>v$ in $H$. If $u'>u$ then by definition $lex(u',v')>lex(u,v)$ (or $lex(\beta) > lex(\alpha)$). But if $u'\leq u$ in $G$ and $v'>v$ in $H$ then since $(u,v)$ is adjacent to $(u',v')$ we obtain $u=u'$ and $v'>v$. It follows that in this case too $lex(\beta) > lex(\alpha)$. It implies that $\{\alpha\}\cup N$ is a descent in $(\mathcal{C}, lex)$, a contradiction. Hence $\mathcal{C}$ is descent-free in $(G\Box H, \tau, \mathcal{C})$. By applying Theorem \ref{descent} for $(\tau, \mathcal{C})$ we obtain
$\mathcal{C}=\mathcal{C}'$. In particular $FF(G\Box H, lex)=FF(G\Box H, \tau)$.
\end{proof}

\noindent We devote the rest of this section to study the First-Fit coloring of $(G\Box H, lex)$. The First-Fit coloring of $(K_m\Box K_n, lex)$ has a significant application in First-Fit coloring of $(G\Box H, lex)$. Hence we begin with an elementary but useful result concerning the Grundy number of $K_m\Box K_n$.

\begin{prop}
$\Gamma(K_m\Box K_n)=\left\{\begin{array}{ll}m+n-1,&m<n\\
2n-2,&m=n.\end{array}\right.$\label{complete}
\end{prop}

\noindent \begin{proof}
Assume that $m<n$. Since $\Delta(K_m\Box K_n)=m+n-2$ then it is enough to obtain a First-Fit coloring of $K_m\Box K_n$ using $m+n-1$ colors. Let $(i,j)$ be any typical vertex of the graph. For each $i,j$ with $1\leq i \leq m$ and $1\leq j \leq n-1$, let the color of $(i,j)$ be $i+j-1~(mod~n-1)$ (replace the color 0 by $n-1$). The result is a pre-coloring of the first $n-1$ columns of $K_m\Box K_n$ using $n-1$ colors such that any vertex of color say $r$ has a neighbor with color $s$ for any $s$ with $s <r$. Now color greedily the vertices of the last column from up to down. We obtain a First-Fit coloring using exactly $m+n-1$ colors. Assume now that $m=n$. In this case using the previous part, it is enough to prove that no First-Fit coloring of $K_n\Box K_n$ uses $2n-1$ colors. Assume on the contrary that ${\mathcal{C}}$ is a First-Fit coloring of the graph using $2n-1$ color and let $v$ be any vertex of color $2n-1$ in ${\mathcal{C}}$. Let $i$ and $j$ be the row and column of $v$ in $K_n\Box K_n$, respectively. Then all other vertices in row $i$ or $j$ have distinct colors from $\{1, \ldots, 2n-2\}$. Assume without loss of generality that the only vertex of color 1 is placed in row $i$. The $n-1$ vertices in column $j$ needs a neighbor with color 1. These vertices of color 1 needs $n-1$ distinct columns. But there are only $n-2$ available columns, a contradiction.
\end{proof}

\begin{thm}
Let $(G, \sigma)$ and $(H, \sigma')$ be two pairs of ordered graphs. Let also $FF(G, \sigma)=p$ and $FF(H, \sigma')=q$. Then $$FF(G\Box H, lex)=FF(K_p\Box K_q,lex).$$\label{completecomplete}
\end{thm}

\noindent \begin{proof} We prove the theorem by induction on $|G|+|H|$. The minimum possible value for $|G|+|H|$ is 2. The assertion obviously holds when $|G|+|H|=2$. Assume (induction hypothesis) that the assertion holds for all graphs $G'$ and $H'$ such that $|G'|+|H'|<|G|+|H|$. Consider now $(G,\sigma)$ and $(H,\sigma')$. If both of these graphs are complete graph then the assertion trivially holds. Assume without generality that $G$ is not complete. Let $C_1, \ldots, C_p$ be the color classes obtained by the First-Fit coloring of $(G, \sigma)$. Since $G$ is not complete, then at least one of the color classes has more than one vertex. Let $C_k$ be a color class such that for any $i<k$, $|C_i|=1$ and $|C_k|\geq 2$. Note that for any $i <k$ and any $v\in C_k$ there exists a vertex $u\in C_i$ such that $v$ is adjacent to $u$ and $\sigma(u)<\sigma(v)$. Since for any $i<k$, $|C_i|=1$, we obtain the following fact.

\noindent {\bf Fact.} For any $u\in C_1 \cup \cdots \cup C_{k-1}$ and any $v\in C_k$ we have $\sigma(u)<\sigma(v)$.

\noindent There are two possibilities concerning the classes $C_k, \ldots, C_p$ and the ordering $\sigma$.

\noindent {\bf Case 1.} For any $v_1, v_2\in G$, if $v_1\in C_k$ and $v_2 \in C_{k+1}\cup \cdots \cup C_p$, then $\sigma(v_1)<\sigma(v_2)$.

\noindent In this case let $G'$ be the graph obtained from $G$ by identifying all vertices in $C_k$ into one vertex, say $w$. Let $\tau$ be the ordering of the vertices in $G'$ obtained by the restriction of $\sigma$ on $G'$. In fact by the above fact and the conditions of Case 1, for any $v\in C_{k+1}\cup \cdots \cup C_p$ (resp. $v \in C_1\cup \cdots \cup C_{k-1}$) we have $\tau(w)<\tau(v)$ (resp. $\tau(w) > \tau(v)$). It is clear that $FF(G', \tau)=p$ and $|G'|<|G|$. Also, by applying the induction hypothesis for $G'\Box H$ we have $FF(K_p\Box K_q,lex)=FF(G'\Box H, lex)$. We show that $FF(G'\Box H, lex)=FF(G\Box H,lex)$. For any $v\in V(G)$, the subgraph of $G\Box H$ induced by $\{(v,u): u\in V(H)\}$ is isomorphic to $H$. Denote this subgraph by $H(v)$. Let $v$ and $v'$ be two arbitrary vertices of $C_k$. Because of the conditions in Case 1 and that $v$ and $v'$ are not adjacent, we obtain that the colorings of $H(v)$ and $H(v')$ in the First-Fit coloring of $(G\Box H, lex)$ are the same. Now we collapse $C_k$ to obtain $G'\Box H$ and its corresponding First-fit coloring. This in particular shows that $FF(G'\Box H, lex)=FF(G\Box H,lex)$. This completes the proof in this case.

\noindent {\bf Case 2.} There exists $v_1\in C_k$ and $v_2 \in C_{k+1}\cup \ldots \cup C_p$ such that $\sigma(v_2)<\sigma(v_1)$.

\noindent In this case we change the order of $v_1$ in $\sigma$ as follows. Put the vertex $v_1$ before (with respect to $\sigma$) all vertices of $\bigcup_{i=k+1}^p C_i$ and after all vertices in $(\bigcup_{i=1}^k C_i)\setminus \{v_1\}$. Denote the new ordering by $\tau$. The following hold for the ordering $\tau$. For any $v\in \bigcup_{i=1}^k C_i$, $\tau(v)<\tau(v_1)$ and for any $u\in \bigcup_{i=k+1}^p C_i$, $\tau(u)>\tau(v_1)$. Also for any $u, v \in V(G)\setminus \{v_1\}$, $\sigma(u)<\sigma (v)$ if and only if $\tau(u)<\tau(v)$. We make the following two claims.

\noindent {\bf Claim 1.} The color classes in the First-Fit coloring of $(G,\tau)$ is the same as the color classes in the First-Fit coloring of $(G,\sigma)$. 

\noindent {\bf Proof of Claim 1.} For simplicity denote the First-Fit coloring of $(G, \sigma)$ and $(G, \tau)$ by ${\mathcal{C}}$ and ${\mathcal{C}}'$, respectively. Then $C_1, \ldots, C_p$ are the color classes in ${\mathcal{C}}$. First, note that any vertex in $C_1$ gets the color $1$ in ${\mathcal{C}}'$. Inductively, any vertex in $C_1 \cup \ldots \cup C_k$ receives the same color in ${\mathcal{C}}'$, since the only vertex whose order is changed is $v_1$ and $v_1\in C_k$. The vertex $v_1$ has at least $k-1$ distinct neighbors of colors $1, \ldots, k-1$ (in the coloring ${\mathcal{C}}$). These neighbors are still before $v_1$ (in the ordering $\tau$) and have the distinct colors $1, \ldots, k-1$ (in the coloring ${\mathcal{C}}'$). It follows that the color of $v_1$ in ${\mathcal{C}}'$ is $k$, as before. Now let $u$ be any arbitrary vertex (including the vertex $v_2$) whose color (in ${\mathcal{C}}$) is $i$, for some $i>k$. We may assume by the induction that for any $1\leq j \leq i-1$, the color (in ${\mathcal{C}}'$) of any vertex in $C_j$ is $j$. The order of $v_1$ (with respect to $\tau$) is before any vertex in $C_{k+1} \cup \ldots \cup C_p$. We conclude that the color of $u$ in ${\mathcal{C}}'$ is $i$. This completes the proof of Claim 1.

\noindent Let $lex'$ be the lexicographic order induced by $(G,\tau)$ and $(H, \sigma')$. Recall that $lex$ is the lexicographic order corresponding to $(G, \sigma)$ and $(H, \sigma')$. We consider $G\Box H$ with two orderings $lex$ and $lex'$. 

\noindent {\bf Claim 2.}
$FF(G\Box H, lex')=FF(G\Box H,lex)$. 

\noindent {\bf Proof of Claim 2.}
Consider the vertex set of $G\Box H$ as a $|V(G)|\times |V(H)|$ array (denoted by ${\mathcal{A}}$) of vertices, where the rows are indexed by the vertices of $G$ and are ordered according to the $\sigma$ from the smallest order at top to the highest order in down of the rows. Also the columns of the array are indexed by $\sigma'$ from the smallest order in left side to the highst order in the right side of the columns. The vertex set of $(G\Box H, lex')$ is the same as this array in which the rows are ordered according to the order $\tau$. It is obvious each vertex in the first row of ${\mathcal{A}}$
gets a same color in the First-Fit colorings of $(G\Box H, lex')$ and $(G\Box H, lex)$. The same is true for the first column of ${\mathcal{A}}$ because of Claim 1. Let $(v,u)$ be any vertex such that $\sigma(v)<\sigma(v_1)$. Let $(v',u')$ be any vertex adjacent to $(v,u)$ such that the order in $lex$ of $(v',u')$ is lower than $(v,u)$. Now, either $v=v'$ and $\sigma'(u')<\sigma'(u)$ or $u=u'$ and $\sigma(v')<\sigma(v)$. Note that in the second case  $\tau(v')<\tau(v)$. It follows that the order in $lex'$ of $(v',u')$ is lower than $(v,u)$. This fact shows that the color of such a vertex $(v,u)$ is identical in both colorings of $G\Box H$. Now consider a vertex of form $(v_1, u)$ and let its color with respect to $lex$ be $t$. Then there are $t-1$ neighbors of $(v_1,u)$ with lower order (in $lex$) and having the colors $1, \ldots, t-1$. Using the above argument we obtain that these neighbor are before (in the ordering $lex'$) the vertex $(v_1,u)$ and hence their colors are identical in the both colorings. It turns out that $(v_1,u)$ gets the color $t$ in the First-Fit coloring of $FF(G\Box H, lex')$. The rest of the vetices are discussed in a similar method and we proceed row by row until are vertices are checked. We obtain that the two First-Fit colorings are identical. The proof of Claim 2 is completed.

\noindent We continue the proof of the theorem. If Case 1 holds for $(G,\tau)$ then using the argument of Case 1 we obtain $FF(G\Box H, lex')=FF(K_p\Box K_q,lex)$. If Case 2 holds for $(G,\tau)$ then we replace $(G,\sigma)$ by $(G,\tau)$ and repeat the above technique for $(G,\tau)$. Let $B(\sigma)$ be the set consisting of the vertices $v\in C_k$ such that there exists a vertex $u\in \bigcup_{i=k+1}^p C_i$ with $\sigma(u)<\sigma(v)$. Note that each time we obtain a new ordering $\tau$ from $\sigma$ then $|B(\sigma)|$ strictly decreases. It turns out that by repeating this technique we eventually obtain an ordering $\tau''$ such that $B(\tau'')$ vanishes. This means that Case 1 holds for $(G,\tau'')$. Let $lex''$ be the lexicographic order corresponding to $\tau''$ and $\sigma'$. We finally obtain $FF(G\Box H, lex)=FF(G\Box H, lex'')=FF(K_p\Box K_q,lex)$, as desired.
\end{proof}

\noindent As we mentioned earlier $\Gamma(G\Box G)$ does not admit any upper bound in terms of $\Gamma(G)$. But for $FF(G\Box G, lex)$ we have a much better result.

\begin{thm}
For any $(G,\sigma)$ and $(H,\sigma')$ we have $$\left\{\begin{array}{ll} FF(G\Box H, lex) \leq \Gamma(G)+\Gamma(H)-1,\\
FF(G\Box G, lex) \leq 2\Gamma(G)-2.\end{array}\right.$$\label{bound}
\end{thm}
\noindent \begin{proof}
Let $FF(G, \sigma)=p$, $FF(H, \sigma')=q$ and $p\leq q$. By Theorem \ref{completecomplete} and Proposition \ref{complete} we have the following two lines of inequalities which yield the result. $$FF(G\Box H, lex)=FF(K_p\Box K_q, lex)\leq \Gamma(K_p\Box K_q) =p+q-1\leq \Gamma(G)+\Gamma(H)-1.$$
$$FF(G\Box G, lex)=FF(K_p\Box K_p, lex)\leq \Gamma(K_p\Box K_p) =2p-2\leq 2\Gamma(G)-2.$$
\end{proof}

\noindent Let $\Bbb{Z}_2=\{0,1\}$ be the only group of size two. Let also $\Bbb{G}_t$ be the direct sum of $t$ copies of $\Bbb{Z}_2$, i.e. $\Bbb{G}_t=\Bbb{Z}_2 \bigoplus \ldots \bigoplus \Bbb{Z}_2$. Consider $\Bbb{G}_t$ as an additive group and denote its elements by $0, 1, \ldots, n-1$, where $n=2^t$. Note that the order of each element of $\Bbb{G}_t$ (other than $0$) is two. Let ${\mathcal{A}}_t$ be the Cayley table of $\Bbb{G}_t$. In fact ${\mathcal{A}}_t$ is obtained as follows. Consider an $n\times n$ array whose rows (from up to down) and columns (from left to right) are indexed by $0, 1, \ldots, n-1$. For any row $i$ and column $j$ the value of ${\mathcal{A}}_t$ in position $(i,j)$ is $i+j$, where $+$ stands for the addition operation of $\Bbb{G}_t$. Let ${\mathcal{C}}_t$ be the $n\times n$ array obtained by adding one to any entry of ${\mathcal{A}}_t$. The entry set of ${\mathcal{C}}_t$ is $\{1, \ldots, n\}$. The array ${\mathcal{C}}_3$ is depicted in Figure \ref{cayley}.
\begin{figure}[ht]
	$$\begin{tabular}{|cccccccc|}
	\hline&&&&&&&\\[-.7eM]
	1&2&3&4&5&6&7&8\\[.2eM]
	2&1&4&3&6&5&8&7\\[.2eM]
	3&4&1&2&7&8&5&6\\[.2eM]
	4&3&2&1&8&7&6&5\\[.2eM]
	5&6&7&8&1&2&3&4\\[.2eM]
	6&5&8&7&2&1&4&3\\[.2eM]
	7&8&5&6&3&4&1&2\\[.2eM]
	8&7&6&5&4&3&2&1\\[.2eM]
	\hline
	\end{tabular}$$
	\caption{The array ${\mathcal{C}}_3$}\label{cayley}
\end{figure}

\begin{thm}
Let $(G,\sigma)$ be any ordered graph. Then $$FF(G\Box G, lex)=2^{\lceil \log FF(G, \sigma) \rceil}.$$\label{power}
\end{thm}
\noindent \begin{proof} We first show that $FF(K_n\Box K_n, lex)=2^t$, where $t$ is such that $2^{t-1} <n\leq 2^t$. Note that $t= \lceil \log n \rceil$. Assume first that $n=2^t$. It is easy to check by hand that the $n\times n$ array obtained by the First-Fit coloring of $(K_n\Box K_n, lex)$ is the same as the array ${\mathcal{C}}_t$, where the entries belong to $\{1, \ldots, n\}$. Hence $FF(K_n\Box K_n, lex)=2^t$ in this case. Now let $t$ be such that $2^{t-1} <n < 2^t$. Let ${\mathcal{T}}$ be the array obtained from the first $n$ rows and $n$ columns of ${\mathcal{C}}_t$. Since $2^{t-1} < n < 2^t$ then exactly the entries $1, 2, \ldots, 2^t$ appear in ${\mathcal{T}}$. From the other side, the First-Fit coloring of $(K_n\Box K_n, lex)$, as an $n\times n$ array, is the same as the array ${\mathcal{T}}$. It follows that $FF(K_n\Box K_n, lex)=2^t$, where $t= \lceil \log n \rceil$.

\noindent Now we consider the general ordered graph $(G,\sigma)$. Set for simplicity $n=FF(G, \sigma)$. By Theorem \ref{completecomplete} and the above result for $K_n\Box K_n$ we have the following relations which complete the proof. $$FF(G\Box G, lex)=FF(K_n\Box K_n, lex)=2^{\lceil \log n \rceil}.$$
\end{proof}

\section{Greedy defining sets in $G\Box H$}

\noindent The topic of defining sets is a well-known area of combinatorics and appears in graph colorings, Latin squares, combinatorial designs, etc. There are many papers concerning defining sets. We refer the reader to \cite{HMTZ} and the survey paper \cite{DMRS}. The greedy defining sets in graphs were first defined in \cite{Z1}. The greedy defining sets of Latin squares was studied in \cite{Z3} and then in \cite{vR, Z4}. In this paper we consider this concept for Cartesian product of graphs. The previous definition of greedy defining sets was given for the minimum vertex coloring of an ordered graph $G$ with $\chi(G)$ colors. In this paper we consider greedy defining sets for proper vertex colorings of $G\Box H$ using $k$ colors, where $k\geq \chi(G\Box H)$. We need to define a general notation. Let $(G,\sigma)$ be an ordered graph and $S$ a subset of vertices in $G$. Let also $\mathcal{C}(S)$ be a pre-coloring of the vertices of $S$. By the First-Fit coloring of $(G,\sigma)$ subject to $\mathcal{C}(S)$, we mean the First-Fit coloring of $(G,\sigma)$ such that the colors of the vertices of $S$ is fixed and the algorithm skips the vertices of $S$ while scanning the vertices of $G$.

\begin{defin}
Let $\mathcal{C}$ be a proper vertex coloring of $(G\Box H, \tau)$ using $k$ colors $1, \ldots, k$. Let $S$ be a subset of vertices in $G\Box H$ and $\mathcal{C}(S)$ be the pre-coloring of $S$ obtained by the restriction of $\mathcal{C}$ to $S$. Then $S$ is called a $k$-greedy defining set of $(G\Box H, \tau, \mathcal{C})$ (or simply $k$-GDS) if the First-Fit coloring of $(G\Box H, \tau)$ subject to $\mathcal{C}(S)$ is the same as the coloring $\mathcal{C}$.
\end{defin}

\noindent We have the following theorem concerning the relationship between greedy defining sets and descents. Proof of the following theorem is similar to the proof of similar results in \cite{Z1} and \cite{Z4}. We omit its proof.

\begin{thm}
Let $\mathcal{C}$ be a proper vertex coloring of $(G\Box H, \tau)$ using $k$ colors. Let $S$ be any subset of vertices such that $S$ intersects any descent in $(G\Box H, \tau, \mathcal{C})$. Then the set $S$ with its coloring from $\mathcal{C}$ is a $k$-GDS of $\mathcal{C}$.\label{descent=graph}
\end{thm}

\noindent We apply Theorem \ref{descent=graph} for quasi-lexicographic orderings of $G\Box H$. It implies that study of greedy defining sets with respect to quasi-lexicographic orders is reduced to study of greedy defining sets in lexicographic order.

\begin{thm}
Let $\tau$ be any quasi-lexicographic ordering and $\mathcal{C}$ any proper vertex coloring of $G\Box H$. Then a subset $S$ is a greedy defining set for $(G\Box H, lex, \mathcal{C})$ if and only if it is a greedy defining set for $(G\Box H, \tau, \mathcal{C})$.\label{quasi-gds}
\end{thm}

\noindent \begin{proof}
By the proof of Theorem \ref{descent-free}, The collection of descents in $(G\Box H, lex, \mathcal{C})$ and $(G\Box H, \tau, \mathcal{C})$ are the same. The assertion follows by Theorem \ref{descent=graph}.
\end{proof}

\noindent In the rest of this section we assume that the vertices (resp. entries) of $G\Box H$ (resp. Latin rectangles) are ordered lexicographically. The next theorem shows the application of greedy defining sets of Latin rectangles in greedy defining sets of $G\Box H$. Let $p$ and $q$ are positive integers and $p\leq q$. Recall that a $p\times q$ Latin rectangle is a $p\times q$ array with entries $1, 2, \ldots, q$ such that no entry is repeated in each row and column of $R$.

\begin{thm}
Let $(G,\sigma)$ and $(H,\sigma')$ be any two ordered graphs. Let $C_1, \ldots, C_p$ (resp. $D_1, \ldots, D_q$) be the color classes in the First-Fit coloring of $(G,\sigma)$ (resp. $(H,\sigma')$), where $FF(G,\sigma)=p$ and $FF(H,\sigma')=q$ with $p\leq q$. Let $R$ be any $p\times q$ Latin rectangle whose rows are top-down indexed by $1, \ldots, p$ and columns are left-right indexed by $1, \ldots, q$. Let also $S$ be a greedy defining set for $R$. Then the following set is a greedy defining set for $G\Box H$ using $q$ colors.
$$\bigcup_{(i,j):(i,j)\in S}~C_i\times D_j.$$\label{latin-to-graph-gds}
\end{thm}

\noindent \begin{proof}
Consider the following proper coloring $\mathcal{C}$ for $G\Box H$. For any $i$ and $j$ with $1\leq i \leq p$ and $1\leq j \leq q$, let $e(i,j)$ be the entry of $R$ in the position $(i,j)$. Assign the entry $e(i,j)$ to all vertices in $C_i\times D_j$ as their color. Assume that a typical descent in $R$ has entries in the positions $(i_1,j_1)$, $(i_1, j_2)$ and $(i_2, j_1)$, where $i_1<i_2$ and $j_1<j_2$. We have $e(i_1, j_2)=e(i_2, j_1)<e(i_1, j_1)$. We note by our definition of $\mathcal{C}$ that any vertex $v$ from $C_{i_1}\times D_{j_1}$ together with its all neighbors of color $e(i_1, j_2)$ form a descent in $\mathcal{C}$. Conversely, any descent of $\mathcal{C}$ is obtained by this method from a descent in $R$. Since $S$ is a GDS in R, then it intersects any descent of $R$. The subset of vertices corresponding to $S$ is $\bigcup_{(i,j):(i,j)\in S}~C_i\times D_j$. It follows that the latter subset intersects any descent in $\mathcal{C}$. The assertion follows using Theorem \ref{descent=graph}.
\end{proof}

\noindent Consider the set $S$ and $D=\bigcup_{(i,j):(i,j)\in S}~C_i\times D_j$ in Theorem \ref{latin-to-graph-gds}. We say that $D$ is a subset of $V(G\Box H)$ corresponding to the set $S$ of the Latin rectangle $R$. Denote by $R(D)$ the coloring of $D$ obtained by the entries of $R$. Another way to state Theorem \ref{latin-to-graph-gds} is the following.

\begin{thm}
Let $R$ be any $p\times q$ Latin rectangle with $p\leq q$ and $S$ a GSD for $(R,lex)$. Let $(G,\sigma)$ and $(H,\sigma')$ be two graphs with $FF(G,\sigma)=p$ and $FF(H,\sigma')=q$. Let $D$ be the subset of vertices of $G\Box H$ corresponding to the elements of $S$. Then the First-Fit coloring of $(G\Box H, lex)$ subject to $R(D)$ uses $FF(H,\sigma')=q$ colors.\label{latin-to-graph}
\end{thm}

\noindent The following gives more information about descent-free colorings in $G\Box H$.

\begin{cor}
Consider $(G\Box G,lex)$ obtained from an ordered graph $(G,\sigma)$. Let $\mathcal{C}$ be any proper vertex coloring of $G\Box G$ using $\chi(G\Box G)$ colors. If $\mathcal{C}$ is descent-free then for some $k$, $\chi(G)=FF(G\Box G, lex)=2^k$.
\end{cor}

\noindent \begin{proof}
Since there does not exist any descent in $\mathcal{C}$ then the First-Fit coloring of $(G\Box G, lex)$ is the same as the coloring $\mathcal{C}$, where only $\chi(G)$ colors are used. By Theorem \ref{power}, $FF(G\Box G, lex)=2^{\lceil \log FF(G, \sigma) \rceil}=\chi(G)\leq FF(G, \sigma)$. Therefore $\log FF(G, \sigma)$ is integer. Hence $FF(G, \sigma)$ is a power of two.
\end{proof}

\noindent In Theorem \ref{upperbound} we obtain an upper bound for the size of greedy defining sets in $G\Box H$. For this purpose we need to obtain an upper bound for the size of greedy defining sets in Latin rectangles. Before we state the next theorem we need to introduce an object associated to any Latin rectangle. Let $R$ be any Latin rectangle of size $m\times n$ on the entry set $\{1, 2, \ldots, n\}$. Let $i\in \{1, 2, \ldots, n\}$ be any arbitrary and fixed entry. There are $m$ entries equal to $i$ in $R$. First, a graph denoted by $G[i]$ on these $m$ entries is defined in the following form. Two entries $e_1$ and $e_2$ (which both are the same as $i$ but in different rows and columns) are adjacent if and only if with an additional entry they form a descent in $R$. The disjoint union $\bigcup_{i=1}^n G[i]$ forms a graph on $mn$ vertices which we denote by $G(R)$. In the following we need the following extremal result of P\'{a}l Tur\'{a}n. Let $G$ be any graph on $m$ vertices without any clique of size $c$. Then $G$ has at most $(c-2)m^2/(2c-2)$ edges.

\begin{thm}
Let $R$ be any Latin rectangle of size $m\times n$. Then $(R, lex)$ contains a GDS of size at most $$nm-n+m-1- \frac{m\log (4m-4)}{4}.$$\label{upperboundrec}
\end{thm}

\noindent \begin{proof}
The proof is based on the fact that any vertex cover for $G(R)$ is a GDS for $R$. The size of minimum vertex cover equals $mn$ minus the independence number of $G(R)$. Hence we obtain an upper bound for the independence number of $G(R)$. The number of edges of $G[i]$ is maximized when the $m$ entries of $i$ lie in the northeast-southwest diagonal of $R$ and the maximum possible number of entries greater than $i$ are placed in the top of this diagonal. It turns out that for $i\geq n-m+2$, $G[i]$ has at most ${m \choose 2} - {i-n+m-1 \choose 2}$ edges. Also for $i\leq n-m+1$, $G[i]$ has at most $m(m-1)/2$ edges. Let $f(i)$ be the maximum number of independent vertices in $G[i]$. For any $i\leq n-m+1$ we have $f(i)\geq 1$. But for any $i\geq n-m+2$, since the complement graph of $G[i]$ has not any clique of size $f(i)+1$ then by Tur\'{a}n's result we obtain that $G(i)$ has at most $(f(i)-1)m^2/2f(i)$ edges. Also since $G[i]$ has at most ${m \choose 2} - {i-n+m-1 \choose 2}$ edges then the complement of $G[i]$ has at least ${i-n+m-1 \choose 2}$ edges. It follows that for $i\geq n-m+2$, $f(i)\geq m^2/[m^2-(i-n+m)(i-n+m-1)]$. Using the substitution $i=n-t$ we obtain the following

$$\sum_{i=1}^n f(i) \geq (n-m+1) + \sum_{t=0}^{m-2} \frac{m^2}{m^2-(m-t)(m-t-1)}$$
$$\hspace{-2cm}\geq (n-m+1) + \sum_{t=0}^{m-2} \frac{m}{2t+1}$$
$$\hspace{-1cm}\geq (n-m+1) + \frac{m\log (4m-4)}{4}.$$

\noindent It turns out that $G(R)$ has at least $(n-m+1) + \frac{m\log (4m-4)}{4}$ independent vertices. Therefore $G(R)$ has a vertex cover of size no more than $nm-n+m-1- \frac{m\log (4m-4)}{4}$. This completes the proof.
\end{proof}

\noindent Using this upper bound for Latin rectangles we obtain a bound for $(G\Box H, lex)$.

\begin{thm}
Let $FF(G,\sigma)=p$, $FF(H,\sigma')=q$ and $p\leq q$. Then $(G\Box H, lex)$ has a $q$-GDS of size at most $\alpha(G)\alpha(H)[pq-q+p-1 - (p\log (4p-4))/4]$.\label{upperbound}
\end{thm}

\noindent \begin{proof}
Let $R$ be a $p\times q$ Latin rectangle. By Theorem \ref{upperboundrec}, there exists a $q$-GDS for $R$ with no more than $[pq-q+p-1 - (p\log (4p-4))/4]$ entries. Note that any color class in $G$ (resp. $H$) has at most $\alpha(G)$ (resp. $\alpha(H)$) vertices. By Theorem \ref{latin-to-graph-gds} we obtain a $q$-GDS for $(G\Box H, lex)$ with no more than $\alpha(G)\alpha(H)[pq-q+p-1 - (p\log (4p-4))/4]$ elements.
\end{proof}

\noindent Let $L$ be any Latin square of size $n\times n$ on the entry set $\{1, 2, \ldots, n\}$. Let also $p$ be any positive integer. By $L+p$ we mean the $n\times n$ Latin square on the entry set $\{p+1, p+2, \ldots, p+n\}$ obtained from $L$ by adding $p$ to each entry of $L$.

\begin{figure}[h]
{\footnotesize
$$\begin{tabular}{|cc|c|c|cc|}
\hline
&&&&&\\
&&&&&\\[0.3eM]
&$L_{k-2}+2^{k-1}+2^{k-2}$&$L_{k-2}+2^{k-1}$&$L_{k-2}+2^{k-2}$&$L_{k-2}$&\\[1.5eM]
&&&&&\\[0.3eM]
\hline
&&&&&\\[1.5eM]
&$L_{k-2}+2^{k-1}$&$L_{k-2}+2^{k-1}+2^{k-2}$&$L_{k-2}$&$L_{k-2}+2^{k-2}$&\\[1.5eM]
&&&&&\\[0.5eM]
\hline
&&&&&\\[1.5eM]
&$L_{k-2}+2^{k-2}$&$L_{k-2}$&$L_{k-2}+2^{k-1}+2^{k-2}$&$L_{k-2}+2^{k-1}$&\\[1.5eM]
&&&&&\\[0.5eM]
\hline
&&&&&\\[1.5eM]
&$L_{k-2}$&$L_{k-2}+2^{k-2}$&$L_{k-2}+2^{k-1}$&$L_{k-2}+2^{k-1}+2^{k-2}$&\\[1.5eM]
&&&&&\\[0.0eM]
\hline
\end{tabular}$$
}
\caption{Proof of Theorem \ref{latin}: Decomposition of $L_k$ into subsquares}\label{fig-latin}
\end{figure}

\begin{thm}
Let $(L_k, lex)$ be the following tensor product, where the number of copies is $k$.
{\footnotesize
\begin{center}
$L_k$~= \begin{tabular}{|c|c|}\hline&\\[-.8eM]
$2$&$1$\\[.4eM]
\hline&\\[-.8eM]
$1$&$2$\\[.4eM]
\hline
\end{tabular}~
$\otimes$~
\begin{tabular}{|c|c|}\hline&\\[-.8eM]
$2$&$1$\\[.4eM]
\hline&\\[-.8eM]
$1$&$2$\\[.4eM]
\hline
\end{tabular}~ $\otimes$ ~$\cdots$~ $\otimes$~
\begin{tabular}{|c|c|}\hline&\\[-.8eM]
$2$&$1$\\[.4eM]
\hline&\\[-.8eM]
$1$&$2$\\[.4eM]
\hline
\end{tabular}~.
\end{center}
}
\noindent Then $L_k$ contains a greedy defining set of cardinality $n^2-{\Omega}(n^{1.673})$, where $n=2^k$.\label{latin}
\end{thm}

\noindent \begin{proof}
Note that $L_k$ is a symmetric array with respect to the two main diagonals of the array. If we divide $L_k$ into four equal subsquares then the north-east and south-west subsquares are equal to $L_{k-1}$. And the north-west and south-east subsquares are equal to $L_{k-1}+2^{k-1}$. We conclude by the induction on $k$ that $L_k$ is decomposed into the sixteen subsquares of size $2^{k-2}\times 2^{k-2}$, as displayed in Figure \ref{fig-latin}.
\noindent Note that the minimum greedy defining number of these sixteen subsquares are all equal. In the following we obtain a greedy defining set, denoted by $D_k$ for $L_k$. For $k=0,1$, $D_0$ is an empty set and $D_1$ consists of a single entry. For $k=2$, $D_2$ is displayed in Figure \ref{small-fig}. Assume that we have obtained $D_{k-1}$ and $D_{k-2}$, for some $k\geq 3$. The greedy defining set $D_k$ for $L_k$ is obtained as follows. Consider the 16 subsquares of $L_k$ as depicted in Figure \ref{fig-latin} and correspond these 16 subsquares with the 16 entries of $L_2$ as illustrated in Figure \ref{small-fig}. Let $S$ be any typical subsquare of $L_k$. If $S$ corresponds with an entry of $L_2$ which belongs to $D_2$ (except the entry 4 of $D_2$ in position $(3,3)$) then we put all entries of $S$ in $D_k$. The total number of these entries is $2^{2k-2}+2^{2k-4}$. There are now eleven subsquares of $L_k$ which have not yet been considered. The four subsquares of these eleven subsquares form the south-east $2^{k-1}\times 2^{k-1}$ subsquare of $L_k$. This subsquare is the same as $L_{k-1}+2^{k-1}$. We put those entries of $L_{k-1}+2^{k-1}$ which correspond to the entries of $D_{k-1}$. The number of these entries is $|D_{k-1}|$. The remaining seven subsquares in $L_k$ are either $K_{k-2}$, or $L_{k-2}+2^{k-2}$, or $L_{k-2}+2^{k-1}$ and or $L_{k-2}+2^{k-2}+2^{k-1}$. We pick from all of these seven subsquares those entries which correspond to the entries of $L_{k-2}$ and put them in $D_k$. Note that the resulting set $D_k$ is a GDS for $(L_k, lex)$.
Set $d_k=|D_k|$. We have $d_k=2^{2k-2}+2^{2k-4}+d_{k-1}+7d_{k-2}$. Let $f(x)=\sum_{k=0}^{\infty} d_kx^k$. Using the recursive relation we obtain
$$f(x)=\frac{1}{1-4x}-\frac{1}{1-x-7x^2}-\frac{2x}{1-x-7x^2}.$$
Let $\alpha=(1-\sqrt{29})/2$ and $\beta=(1+\sqrt{29})/2$ so that $1-x-7x^2=(1-\alpha x)(1-\beta x)$.
Set $n=2^k$. We obtain $$d_k=n^2-n^{\log \beta}(\frac{\sqrt{29}+5}{2\sqrt{29}})-(-1)^k n^{\log (-\alpha)}(\frac{\sqrt{29}-5}{2\sqrt{29}}).$$
Since $\log \beta$ is approximately 1.6735 then $d_k=n^2-{\Omega}(n^{1.673}).$
\begin{figure}
{\footnotesize
$$\begin{tabular}{|cccc|}
\hline
&&&\\[-0.6eM]
4&$\fbox{3}$&$\fbox{2}$&1\\[0.2eM]
&&&\\[0.2eM]
3&4&$\fbox{1}$&$\fbox{2}$\\[0.2eM]
&&&\\[0.2eM]
2&1&$\fbox{4}$&3\\[0.2eM]
&&&\\[0.2eM]
$\fbox{1}$&2&3&4\\[0.2eM]
\hline
\end{tabular}$$
}
\caption{$L_2$ with a GDS of size 6}\label{small-fig}
\end{figure}
\end{proof}

\noindent Let $m_k$ be the size of minimum GDS in $L_k$. It is easily seen that $m_k\geq 4m_{k-1}$. Also $m_2=6$. We obtain the following result.

\begin{thm}
Let $n=2^k$. Then any GDS in $L_k$ needs at least $6n^2/16$ entries.
\end{thm}

\section{Questions for further researches}

\noindent In this section we propose some questions for further researches. All greedy defining sets in this section are considered for lexicographic order. In the previous section we tried to obtain the best possible upper bound for the minimum greedy defining sets in all Latin squares. Because $L_k$ (see Theorem \ref{latin}) has maximum number of descents in all known families of Latin squares, we guess that probably the minimum GDS of $L_k$ has the maximum value among all Latin squares of size $n=2^k$. But Theorem \ref{latin} is our best result and we could not obtain a GDS for $L_k$ with cardinality at most $\lambda n^2$, for some constant $\lambda<1$. Hence we propose the following question.

\begin{prob}
Does there exist a constant $\lambda<1$ such that any Latin square of size $n$ has a greedy defining set of cardinality at most
$\lambda n^2$.
\end{prob}

\noindent Let $g_n$ be the cardinality of smallest greedy defining set among all $n\times n$ Latin squares. The following conjecture from \cite{Z3} is still unsolved. Although it was proved for some infinite sequences of natural numbers.

\begin{conj}
$$g_n={\mathcal{O}}(n).$$
\end{conj}

\noindent The next question concerns $L_k$ in Theorem \ref{latin}.

\begin{prob}
Determine the greedy defining number of $(L_k, lex)$.
\end{prob}

\noindent We finally propose the following complexity problem from \cite{Z3}. We conjecture now that the answer is affirmative.

\begin{prob}
Given any Latin square $(L, lex)$, is to determine the minimum cardinality of GDS in $L$, an $\mathcal{N}\mathcal{P}$-complete problem?
\end{prob}

\section{Acknowledgment}

\noindent The author thanks the anonymous referee for kind reviewing of the paper.


\begin{thebibliography}{1}

\bibitem{AHL}
M. Ast\'{e}, F. Havet, C. Linhares Sales, Grundy number and products of
graphs, Discrete Math. 310 (2010), 1482--1490.

\bibitem{AHI}
F, Aurenhammer, J. Hagauer, W. Imrich, Cartesian graph factorization at logarithmic cost per edge, Comput. Complexity 2 (1992), 331--349.

\bibitem{BHLY}
J. Balogh, S. G. Hartke, Q. Liu, G. Yu, On the first-fit chromatic number
of graphs, SIAM J Discrete Math. 22 (2008), 887--900.

\bibitem{CGHLM}
V. Campos, A. Gy\'{a}rf\'{a}s, F. Havet, C. Linhares Sales, F. Maffray, New bounds on the Grundy number of products of graphs, J. Graph Theory 71 (2012), 78--88.

\bibitem{DMRS}
D. Donovan, E. S. Mahmoodian, C. Ramsay, A. P. Street, Defining
sets in combinatorics: a survey, Surveys in combinatorics, 2003
(Bangor), 115--174, London Math. Soc. Lecture Note Ser., 307,
Cambridge Univ. Press, Cambridge (2003)

\bibitem{HMTZ}
H. Hajiabolhassan, M. L. Mehrabadi, R. Tusserkani, M. Zaker, A characterization of uniquely vertex colorable graphs using minimal defining sets, Disc. Math. 199 (1999) 233--236.

\bibitem{HJ}
D. G. Hoffman, Jr. P. D. Johnson, Greedy colorings and the Grundy chromatic number of the $n$-cube, Bull. ICA 26 (1999), 49--57.

\bibitem{HKS}
F. Havet, T. Kaiser, M. Stehlik, Grundy number of the Cartesian product
of a tree and a graph, unpublished manuscript.

\bibitem{vR}
J. van Rees, More greedy defining sets in Latin squares, Australas. J. Combin. 44 (2009), 183--198.

\bibitem{Z1}
M. Zaker, Greedy defining sets of graphs, Australas. J. Combin. 23 (2001), 231--235.

\bibitem{Z2}
M. Zaker, Results on the Grundy chromatic number of graphs, Discrete Math. 306 (2006), 3166--3173.

\bibitem{Z3}
M. Zaker, Greedy defining sets of Latin squares, Ars Combin. 89 (2008), 205--222.

\bibitem{Z4}
M. Zaker, More results on greedy defining sets, Ars Combin. 114 (2014), 53--64.
\end{thebibliography}
\end{document}